\documentclass[24pt,reqno]{article}

\usepackage{amsmath}
\usepackage{amssymb}
\usepackage{color}
\usepackage{epsf}

\begin{document}
\makeatletter

\begin{center}
\epsfxsize=10in
\end{center}

\def\endofproofmark{$\Box$}

\def\endofproofmark{$\Box$}

\begin{center}
{\LARGE\bf Some New Inequalities Between Important Means }\vskip
2mm {\LARGE\bf and Applications to Ky Fan - type Inequalities}
\vskip 1cm
\bf{Jamal Rooin}\\
\bf{Mehdi Hassani}\\
\vskip .5cm
Department of Mathematics\\
Institute for Advanced Studies in Basic Sciences\\
P.O.Box: 45195-1159, Zanjan, Iran\\
\vskip 2mm
{\tt rooin@iasbs.ac.ir}\\
{\tt mhassani@iasbs.ac.ir}\\
\end{center}

\date{}
\newtheorem{theo}{Theorem}[section]
\newtheorem{prop}[theo]{Proposition}
\newtheorem{rem}{Remark}
\newtheorem{lemma}[theo]{Lemma}
\newtheorem{cor}[theo]{Corollary}
\newtheorem{problem}{Problem}
\def\frameqed{\framebox(5.2,6.2){}}
\def\deshqed{\dashbox{2.71}(3.5,9.00){}}
\def\ruleqed{\rule{5.25\unitlength}{9.75\unitlength}}
\def\myqed{\rule{8.00\unitlength}{12.00\unitlength}}
\def\qed{\hbox{\hskip 6pt\vrule width 7pt height11pt depth1pt\hskip 3pt}
\bigskip}
\newenvironment{proof}{\trivlist\item[\hskip\labelsep{\bf Proof}.]}{\hfill
 $\frameqed$ \endtrivlist}
\newcommand{\COM}[2]{{#1\choose#2}}

\thispagestyle{empty} \null \addtolength{\textheight}{1cm}

\begin{abstract}
In this paper, mainly using the convexity of the function
$\frac{a^x-b^x}{c^x-d^x}$ and convexity or concavity of the
function $\ln\frac{a^x-b^x}{c^x-d^x}$ on the real line, where
$a>b\geq c>d>0$ are fixed real numbers, we obtain some important
relations between various important means of these numbers. Also,
we apply the obtained results to Ky Fan type inequalities and get
some new refinements.

\end{abstract}

\bigskip
\hrule
\bigskip



\bigskip
\hrule
\bigskip

\section{Introduction and Motivation}
Suppose that $a>b\geq c>d>0$. It is shown in \cite{rooin} that the
function
$$
f(x)=\frac{a^x-b^x}{c^x-d^x}\hspace{15mm}(-\infty<x<+\infty),
$$
is strictly increasing on the real line, moreover
$\lim\limits_{x\rightarrow +\infty}f(x)=+\infty$ and
$\lim\limits_{x\rightarrow -\infty}f(x)=0$. By a simple
calculation, we have
\begin{eqnarray}
f(x)&=&\left\{ \begin{array}{ll}
\frac{a-b}{c-d}\frac{L(c,d)}{L(a,b)} & x=0,\\
\frac{a-b}{c-d} & x=1,\\
\frac{a-b}{c-d}\left(\frac{L_{x-1}(a,b)}
{L_{x-1}(c,d)}\right)^{x-1} & x\neq 0,1,
\end{array}
\right.\nonumber\\
\frac{f'(x)}{f(x)}&=&\frac{1}{x}\ln\frac{I(a^x,b^x)}{I(c^x,d^x)}\hspace{15mm}(x\neq
0),\\
f'(0)&=&\frac{a-b}{c-d}\frac{L(c,d)}{L(a,b)}\ln\frac{G(a,b)}{G(c,d)}.\nonumber
\end{eqnarray}
Note that the notations $L(a,b),~L_{p}(a,b)$ and $G(a,b)$ are
well-known means between $a,b>0$. We recall them in the following
table:
\begin{center}
\begin{tabular}{|c|c|c|}
  \hline
  Name & Notation & Definition \\
  \hline
  arithmetic mean & $A(a,b)$ & $\frac{a+b}{2}$ \\
  \hline
  geometric mean & $G(a,b)$ & $\sqrt{ab}$ \\
  \hline
  harmonic mean & $H(a,b)$ & $\frac{2}{\frac{1}{a}+\frac{1}{b}}$ \\
  \hline
  logarithmic mean & $L(a,b)$ & $\left\{ \begin{array}{ll}
a & a=b \\
\frac{a-b}{\ln a-\ln b} & a\neq b
\end{array}
\right.$ \\
\hline
identric mean & $I(a,b)$ & $\left\{ \begin{array}{ll}
a & a=b \\
\frac{1}{e}\left(\frac{a^a}{b^b}\right)^{\frac{1}{a-b}} & a\neq b
\end{array}
\right.$ \\
\hline
 $p-$logarithmic mean & $L_p(a,b)$ & $\left\{ \begin{array}{ll}
a & a=b \\
\left(\frac{a^{p+1}-b^{p+1}}{(p+1)(a-b)}\right)^{\frac{1}{p}} &
a\neq b
\end{array}
\right.\hspace{15mm}p\neq 0,-1$\\
  \hline
\end{tabular}
\end{center}
\vspace{2mm}
\begin{rem}
{\rm (i)} With above notations, we have the following limit cases:
\begin{eqnarray}
\lim_{p\rightarrow 0}L_p(a,b)=I(a,b),\hspace{15mm}
\lim_{p\rightarrow -1}L_p(a,b)=L(a,b).
\end{eqnarray}
\hspace{21mm}{\rm (ii)} The following inequalities are well-known in
the literature \cite{bul-mit-vas}:
\begin{eqnarray}
H(a,b)\leq G(a,b)\leq L(a,b)\leq I(a,b)\leq A(a,b),
\end{eqnarray}
and equality holds in each inequality if and only if $a=b$.
\end{rem}
In \cite{bbr} it is declared that $f$ is strictly convex on the real
line, which by considering $f(x)\rightarrow 0~(x\rightarrow
-\infty)$, it is a stronger result than being strictly increasing,
and besides, the function
$$
g(x)=\ln\frac{a^x-b^x}{c^x-d^x},
$$
is strictly convex if $ad-bc>0$, and is strictly concave if
$ad-bc<0$. Since, we can write $f$ and $g$ in terms of means, we
can use the convexity or concavity of them in order to get some
relations between different means mentioned in the above
table.\\
In this paper, first we study these functions more closely and get
some interesting inequalities which are contained in the heart of
these functions, and then as applications, using the achieved
results, we sharpen some Ky Fan type inequalities. We refer the
readers who are interested in further results on the $p-$logarithmic
mean and in Ky Fan - type inequalities to \cite{alz1, le-sh1,
le-sh2, le-sh3}.

\section{Study of the Functions $\frac{a^x-b^x}{c^x-d^x}$ and $\ln\frac{a^x-b^x}{c^x-d^x}$}
In this section, we will prove our claims about the functions $f$
and $g$. First $g$:
\begin{theo}
Suppose $a>b\geq c>d>0$, and let
$$g(x)=\ln\frac{a^x-b^x}{c^x-d^x}.$$
Then $g$ is strictly convex if $ad-bc>0$, and is strictly concave
if $ad-bc<0$. If $ad-bc=0$, then $g$ turns out to be a linear
mapping.
\end{theo}
\begin{proof}Suppose that $ad-bc>0$. Since
$$
g(x)=\ln\frac{(\frac{a}{b})^x-1}{(\frac{c}{d})^x-1}+x\ln\frac{b}{d},
$$
it is sufficient to show that if $a>b>1$, then
$\ln\frac{a^x-1}{b^x-1}$ is strictly convex, and since
$\ln\frac{a^x-1}{b^x-1}=\ln\frac{e^{x\ln a}-1}{e^{x\ln b}-1}$, it
is sufficient to show that if $a>b>0$, then
$$
u(x)=\ln\frac{e^{ax}-1}{e^{bx}-1}
$$
is strictly convex. A simple calculation yields that
$$
u''(x)=\frac{b^2 e^{bx}}{(e^{bx}-1)^2}-\frac{a^2
e^{ax}}{(e^{ax}-1)^2}\hspace{15mm}(x\neq 0).
$$
So, for $x\neq 0$, $u''(x)>0$ is equivalent to
$\frac{|\sinh\frac{ax}{2}|}{a}>\frac{|\sinh\frac{bx}{2}|}{b}$, or
$\frac{\sinh\frac{ax}{2}}{\frac{ax}{2}}>\frac{\sinh\frac{bx}{2}}{\frac{bx}{2}}$.
But it is clear that the function $\frac{\sinh x}{x}$ is strictly
decreasing on $(-\infty,0]$ and strictly increasing on
$[0,+\infty)$. This yields our claim in this case.\\
In the case $ad-bc<0$, rewrite $g$ as follows
$$
g(x)=-\ln\frac{(\frac{c}{d})^x-1}{(\frac{a}{b})^x-1}+x\ln\frac{b}{d}.
$$
According to the above argument, the function
$-\ln\frac{(\frac{c}{d})^x-1}{(\frac{a}{b})^x-1}$, and so $g$, is
strictly concave.\\
If $ad-bc=0$, then $g(x)=x\ln\frac{b}{d}$; a straight line through
the origin. This completes the proof.
\end{proof}

Now, consider the function $f$. If $ad-bc=0$, then
$f(x)=(\frac{b}{d})^x$ which is clearly strictly  convex. If
$ad-bc>0$, $g$ is strictly convex, and since the function "exp" is
strictly increasing and convex, $f=\exp(\ln(f))=\exp(g)$ is
strictly convex. But, in the case of $ad-bc<0$, we cannot use the
above method. Therefore, we go to prove the convexity of $f$
independently.
\begin{theo}
Suppose $a>b\geq c>d>0$ and
$$
f(x)=\frac{a^x-b^x}{c^x-d^x}.
$$
Then $f$ is strictly convex on the real line.
\end{theo}
\begin{proof}
Since $\frac{a}{d}>\frac{b}{d}\geq\frac{c}{d}>1$ and
$f(x)=\frac{(\frac{a}{d})^x-(\frac{b}{d})^x}{(\frac{c}{d})^x-1}$,
it is sufficient to consider
$$
f(x)=\frac{a^x-b^x}{c^x-1}\hspace{15mm}(a>b\geq c>1),
$$
and since $f(x)=\frac{e^{x\ln a}-e^{x\ln b}}{e^{x\ln c}-1}$ and
$\ln a>\ln b\geq \ln c>0$, it is sufficient to consider
$$
f(x)=\frac{e^{ax}-e^{bx}}{e^{cx}-1}\hspace{15mm}(a>b\geq c>0).
$$
But, $f(x)=\frac{e^{ax}}{e^{cx}-1}-\frac{e^{bx}}{e^{cx}-1}$, and
so
$$
f''(x)=\frac{e^{ax}\{[(a-c)e^{cx}-a]^2+c^2e^{cx}\}-e^{bx}\{[(b-c)e^{cx}-b]^2+c^2e^{cx}\}}
{(e^{cx}-1)^3}.
$$
Therefore, it is sufficient to prove that for any fixed
$x>0~(x<0)$, the function
$$
h(t)=e^{tx}\{[(t-c)e^{cx}-t]^2+c^2e^{cx}\},
$$
is strictly increasing (decreasing) on $t\geq c$. But,
$$
h'(t)=xe^{tx}\{[(t-c)e^{cx}-t]^2+c^2e^{cx}\}+2(e^{cx}-1)[(t-c)e^{cx}-t]e^{tx}.
$$
Since $e^{tx}>0$, the sign of $h'(t)$ agrees with the sign of
$$
k(t)=\frac{h'(t)}{e^{tx}}=x\{[(t-c)e^{cx}-t]^2+c^2e^{cx}\}+2(e^{cx}-1)[(t-c)e^{cx}-t].
$$
But,
$$
k'(t)=2x(e^{cx}-1)[(t-c)e^{cx}-t]+2(e^{cx}-1)^2,
$$
and
$$
k''(t)=2x(e^{cx}-1)^2.
$$
So, if $x>0$, then $k'(t)$ is strictly increasing, and so for
$t>c$, $k'(t)>k'(c)$. Also, if $x<0$, then $k'(t)$ is strictly
decreasing and so for $t>c$, $k'(t)<k'(c)$. But,
$k'(c)=2(e^{cx}-1)(e^{cx}-1-cx)$, and $e^{cx}-1-cx>0~(x\neq 0)$.
So, if $x>0$, then $k'(c)>0$ and if $x<0$, then $k'(c)<0$. Thus,
when $x>0$, we have $k'(t)>k'(c)>0~(t>c)$ and when $x<0$ we have
$k'(t)<k'(c)<0~(t>c)$. Therefore, if $x>0$, then $k(t)$ is
strictly increasing and if $x<0$, then $k(t)$ is strictly
decreasing on $[c,\infty)$. So, for $t>c$, if $x>0$, then
$k(t)>k(c)$ and if $x<0$, then $k(t)<k(c)$. Now, let
$$
u(x)=k(c)=xc^2+xc^2e^{cx}-2c(e^{cx}-1).
$$
We have
$$
u'(x)=xc^3e^{cx}-c^2e^{cx}+c^2,
$$
and
$$
u''(x)=xc^4e^{cx}.
$$
Thus if $x>0$, then $u''(x)>0$, $u'(x)>u'(0)=0$, and therefore,
$u(x)>u(0)=0$. Also, if $x<0$, then $u''(x)<0$ and
$u'(x)<u'(0)=0$, and therefore $u(x)<u(0)=0$. So for $t>c$, if
$x>0$, then $k(t)>k(c)=u(x)>0$ and if $x<0$, then
$k(t)<k(c)=u(x)<0$. So, for $t>c$, the sign of $h'(t)$ is same as
the sign of $x$. Thus if $x>0$, then the function $h(t)$ is
strictly increasing on $t\geq c$ and if $x<0$, then the function
$h(t)$ is strictly decreasing on $t\geq c$. So, if $x>0$,
$h(a)>h(b)$ and therefore $f''(x)>0$, and besides if $x<0$,
$h(a)<h(b)$ and again $f''(x)>0$. This completes the proof.
\end{proof}

\section{Applications to Special Means}

As we said before, the convexity of $f$ is a strong equipment for
establishing interesting inequalities. For example, if $a>b\geq
c>d>0$, we have the following nontrivial inequality
$$
\frac{\frac{a^b-b^b}{c^b-d^b}-\frac{a^d-b^d}{c^d-d^d}}{b-d}<
\frac{\frac{a^a-b^a}{c^a-d^a}-\frac{a^c-b^c}{c^c-d^c}}{a-c},
$$
since $m_{d,b}<m_{c,a}$, where $m_{\alpha,\beta}$ denotes the
slope of the line segment joining points $(\alpha,f(\alpha))$ and
$(\beta,f(\beta))$. Some other results are given in the next
theorems.
\begin{theo}
Suppose $a>b\geq c>d>0$ and $p,q \neq0,-1$. Then we have the
following inequality
\begin{eqnarray}
\frac{L_{p}^{p}(a,b)}{L_{p}^{p}(c,d)}\geq
\frac{L_{q}^{q}(a,b)}{L_{q}^{q}(c,d)}\left(1+\left(\frac{p-q}{q+1}\right)
\ln\frac{I(a^{q+1},b^{q+1})}{I(c^{q+1},d^{q+1})}\right),
\end{eqnarray}
with equality holding if and only if $p=q$.\\ In particular
\begin{eqnarray}
\exp\left(1-\frac{L(c,d)}{L(a,b)}\right)<\frac{I(a,b)}{I(c,d)}<
\exp\left(\frac{L(a,b)}{L(c,d)}-1\right),
\end{eqnarray}
and for $a>b>0$,
\begin{eqnarray}
\exp\left(1-\frac{b}{L(a,b)}\right)<\frac{I(a,b)}{b}<
\exp\left(\frac{L(a,b)}{b}-1\right).
\end{eqnarray}
\end{theo}
\begin{proof}
Since $f$ is strictly convex, we have
$$
f(p+1)\geq f(q+1)+(p-q)f'(q+1),
$$
with equality holding if and only if $p=q$. Considering this fact and (1) yields (4).\\
Now, if we put $p=-1$ and $q=0$ in the last relation, considering
(1), we have
$$
\frac{L(c,d)}{L(a,b)}>1-\ln\frac{I(a,b)}{I(c,d)},
$$
which yields the left hand side of (5). Similarly, if we put
$p=-1$ and $q=-2$, then
$$
\frac{L(c,d)}{L(a,b)}>\frac{L_{-2}^{-2}(a,b)}{L_{-2}^{-2}(c,d)}
\left(1-\ln\frac{I(a^{-1},b^{-1})}{I(c^{-1},d^{-1})}\right).
$$
But, $L_{-2}^{-2}(a,b)=\frac{1}{ab}$ and
$\frac{1}{d}>\frac{1}{c}\geq \frac{1}{b}>\frac{1}{a}>0$. So,
changing $\frac{1}{d},~\frac{1}{c},~\frac{1}{b}$ and $\frac{1}{a}$
by $a,~b,~c$ and $d$ respectively and considering
$L(a^{-1},b^{-1})=\frac{1}{ab}L(a,b)$, we get the right hand side
of (5).\\
Since
\begin{eqnarray}
\frac{L(a,b)}{b}=\frac{\frac{a}{b}-1}{\ln\frac{a}{b}}~~~{\rm
and}~~~
\ln\frac{I(a,b)}{b}=-1+\frac{\frac{a}{b}}{\frac{a}{b}-1}\ln\frac{a}{b},
\end{eqnarray}
putting $x=\frac{a}{b}$, the inequalities in (6) follow from
$$
x\ln x+\ln x-2x+2>0\hspace{15mm}(x>1),
$$
$$
(x-1)^2-x\ln^2 x>0\hspace{15mm}(x>1),
$$
respectively.
\end{proof}
\begin{theo}
If $a>b\geq c>d>0$, then
\begin{eqnarray}
\frac{L(a,b)}{L(c,d)}>1+\ln\frac{G(a,b)}{G(c,d)}>\frac{2ab}{ab+cd},
\end{eqnarray}
and
\begin{eqnarray}
\frac{L(a,b)}{L(c,d)}>\frac{\ln\frac{G(a,b)}{G(c,d)}}{\ln\frac{I(a,b)}{I(c,d)}}.
\end{eqnarray}
In particular, if $a>b>0$, then
\begin{eqnarray}
\frac{L(a,b)}{b}>
1+\frac{1}{2}\ln\frac{a}{b}>\frac{2a}{a+b}>\frac{\ln\frac{a}{b}}{2\ln\frac{I(a,b)}{b}}
.
\end{eqnarray}
\end{theo}
\begin{proof}
Since $f(x)=\frac{a^x-b^x}{c^x-d^x}$ is strictly convex, we have
$$
f'(0)x+f(0)<f(x)\hspace{15mm}(x\neq 0),
$$
which by setting $x=1$ and using (1), we get the first inequality
in (8). The second inequality in (8) is equivalent to
\begin{eqnarray}
\ln\frac{G(a,b)}{G(c,d)}>\frac{ab-cd}{ab+cd},
\end{eqnarray}
which by putting $x=ab$ and $y=cd$, (11) follows from
\begin{eqnarray}
\frac{1}{2}\ln\frac{x}{y}>\frac{\frac{x}{y}-1}{\frac{x}{y}+1}\hspace{15mm}(x>y>0).
\end{eqnarray}
But, (12) is obtained by the facts that the function
$h(x)=\frac{1}{2}\ln x-\frac{x-1}{x+1}$ is strictly increasing on
$[1,\infty)$ and $\frac{x}{y}$ is greater than 1.\\
The inequality (9) follows from $f'(0)<f'(1)$ and considering (1).\\
Considering (7) and putting $x=\frac{a}{b}$, the inequalities in
(10) follow from left to right from
$$
\ln^2 x+2\ln x-2x+2<0\hspace{15mm}(x>1),
$$
$$
x\ln x+\ln x-2x+2>0\hspace{15mm}(x>1),
$$
$$
(x^2-1)\ln x-4x(1-x+x\ln x)<0\hspace{15mm}(x>1),
$$
respectively. This completes the proof.
\end{proof}
Now, consider the function $g(x)=\ln\frac{a^x-b^x}{c^x-d^x}$. It
is evident from (1), $g'(0)=\ln\frac{G(a,b)}{G(c,d)}$ and for
$x\neq 0$, $g'(x)=\frac{1}{x}\ln\frac{I(a^x,b^x)}{I(c^x,d^x)}$.
\begin{theo}
Suppose $a>b\geq c>d>0$ and $p,q \neq0,-1$. If $ad-bc>0$, then
\begin{eqnarray}
\left(\frac{L_p(a,b)}{L_p(c,d)}\right)^p\geq
\left(\frac{L_q(a,b)}{L_q(c,d)}\right)^q\left(
\frac{I(a^{q+1},b^{q+1})}{I(c^{q+1},d^{q+1})}\right)^{\frac{p-q}{q+1}}.
\end{eqnarray}
If $ad-bc<0$, the inequality reverses. The equality holds if and
only if $~ad-bc=0$ or $p=q$.
\end{theo}
\begin{proof}
If $ad-bc>0$, then by Theorem 2.1, $g$ is strictly convex and so
considering the tangent line at $x=q+1$, we have
$$
g(p+1)\geq g(q+1)+(p-q)g'(q+1),
$$
with equality holding if and only if $p=q$. Now, considering (1),
we get (13) with equality if and only if $p=q$.\\
If $ad-bc<0$, then $g$ is strictly concave and the argument is
similar.\\
If $ad-bc=0$, then $g(x)=x\ln\frac{b}{d}$ is linear, and so
equality always holds in (13).
\end{proof}
In the cases $p,q=0,-1$, we conclude the following nice result:
\begin{theo}
Suppose $a\geq b\geq c\geq d>0$. If $ad-bc>0$, then
\begin{eqnarray}
\frac{H(a,b)}{H(c,d)}<\frac{G(a,b)}{G(c,d)}<\frac{L(a,b)}{L(c,d)}<
\frac{I(a,b)}{I(c,d)}<\frac{A(a,b)}{A(c,d)}.
\end{eqnarray}
If $ad-bc<0$, all inequalities reverse, and if $ad-bc=0$, all
inequalities turn out to be equalities.
\end{theo}
\begin{proof}
\textbf{Case I.} $ad-bc>0$. It is divided into two branches; $a>b
\geq c>d$ and $a>b\geq c=d$.\\
If $a>b\geq c>d$, writing the first inequality in (14) in terms of
$\frac{a}{b}$ and $\frac{c}{d}$, it follows from the fact that the
function $x+\frac{1}{x}$ is strictly increasing on $[1,\infty)$.\\
The second one follows from the fact that the slope of the line
segment between $(-1,g(-1))$ and $(0,g(0))$ is strictly less than
the
slope of the line segment between $(0,g(0))$ and $(1,g(1))$.\\
The third one follows from the fact that the point $(0,g(0))$ is
strictly above the tangent line to the graph of $g$ at $x=1$.\\
Writing the last inequality in (14) in terms of $\frac{a}{b}$ and
$\frac{c}{d}$, and considering (7), it follows from the fact that
the function $\frac{x\ln x}{x-1}-\ln(x+1)$ is strictly decreasing
on $[1,\infty)$.\\
If $a>b\geq c=d$, all denominators in (14) are all equal to $c$,
and so (14) follows from Remark 1, (ii).\\
\textbf{Case II.} $ad-bc<0$. We have $a>b\geq c>d$ or $a=b\geq
c>d$,
and the result follows similarly by using the strict concavity of $g$.\\
\textbf{Case III.} $ad-bc=0$. It turns to two branches; $a>b \geq
c>d$ and $a=b=c=d$. To prove the first case, proceed as the case I
and use the linearity of $g$. In the second case, all the
fractions in (14) are equal to 1.
\end{proof}
Now, we give a nice example concerning some numerical
sequences.\\\\
\textbf{Example.} For every $n\in\mathbb{N}$, we have
\begin{eqnarray}
\frac{n+2}{n+1}<1+\ln\sqrt{\frac{n+2}{n}}<\frac{\ln(1+\frac{1}{n})}{\ln(1+\frac{1}{n+1})},
\end{eqnarray}
and
\begin{eqnarray}
\frac{\ln\sqrt{\frac{n+2}{n}}}{\ln\frac{(n+2)(1+\frac{1}{n+1})^{n+1}}{(n+1)(1+\frac{1}{n})^n}}
<\frac{\ln(1+\frac{1}{n})}{\ln(1+\frac{1}{n+1})}.
\end{eqnarray}
Also, we have
\begin{eqnarray}
\frac{2n+3}{2n+1}<\frac{n+2}{n+1}\frac{(1+\frac{1}{n+1})^{n+1}}{(1+\frac{1}{n})^{n}}
<\frac{\ln(1+\frac{1}{n})}{\ln(1+\frac{1}{n+1})}<\sqrt{\frac{n+2}{n}}<\frac{(n+2)(2n+1)}{n(2n+3)}.
\end{eqnarray}
These are obtained from (8), (9) and Theorem 3.4, by putting
$a=n+2$, $b=c=n+1$ and $d=n$ with considering $ad-bc<0$.

\section{Applications to Ky Fan Type Inequalities}

Throughout this section, given $n$ arbitrary nonnegative real
numbers $x_1, \cdots, x_n$ belonging to $(0,\frac{1}{2}]$, we
denote by $A_n$ and $G_n$, the unweighted arithmetic and geometric
means of $x_1, \cdots, x_n$ respectively, i.e.
$$
A_n=\frac{1}{n}\sum_{i=1}^n x_i,\hspace{10mm}G_n=\prod_{i=1}^n
x_i^{1/n},
$$
and by $A'_n$ and $G'_n$, the unweighted arithmetic and geometric
means of $1-x_1, \cdots, 1-x_n$ respectively, i.e.
$$
A'_n=\frac{1}{n}\sum_{i=1}^n
(1-x_i),\hspace{10mm}G'_n=\prod_{i=1}^n (1-x_i)^{1/n}.
$$
With the above notations, the Ky Fan's inequality \cite{bb} asserts
that:
\begin{eqnarray}
\frac{A'_n}{G'_n}\leq \frac{A_n}{G_n},
\end{eqnarray}
{\it with equality holding if and only if $x_1=\cdots=x_n$.}\\
In 1988, H. Alzer \cite{alz2} obtained an additive analogue of Ky
Fan's inequality as follows:
\begin{eqnarray}
A'_n-G'_n\leq A_n-G_n,
\end{eqnarray}
{\it with equality holding if and only if $x_1=\cdots=x_n$.}\\
Also, in 1995, J.E. Pe\v{c}ari\'{c} and H. Alzer \cite{pec-alz},
using the Dinghas Identity \cite{ding}, proved that:
\begin{eqnarray}
{A_n}^n-{G_n}^n\leq {A'_n}^n-{G'_n}^n,
\end{eqnarray}
\textit{in which if $n=1,2$, equality always holds in {\rm (20)},
and if $n\geq 3$, the equality is valid if and only
if $x_1=\cdots =x_n$.}\\
Now, we use the obtained results in Theorems 3.2 and 3.4 and find
some refinements and inverses of Ky Fan's inequality (18) and its
additive analogues (19) and (20). Perhaps, the most interesting
results are:
\begin{eqnarray}
\left(\frac{A'_n}{G'_n}\right)^{A'_n+G'_n}\leq
\left(\frac{A_n}{G_n}\right)^{A_n+G_n},
\end{eqnarray}
\begin{eqnarray}
\left(\frac{A'_n}{G'_n}\right)^{A_n-G_n}\leq
\left(\frac{A_n}{G_n}\right)^{A'_n-G'_n},
\end{eqnarray}
{\it with equality holding if and only if $x_1=\cdots=x_n$.}\\
These are easily obtained from the first inequality in (29) which
need a great labor to handle them directly.
\begin{theo} Suppose $x_1, \cdots,
x_n\in(0,\frac{1}{2}]$ not all equal. Then
\begin{eqnarray}
\frac{A'_n}{G'_n}<
\left(\frac{A_n}{G_n}\right)^{\frac{{A'_n}-{G'_n}}{A_n-G_n}-\frac{\ln\frac{A'_n}{G'_n}}
{\ln\frac{A_n}{G_n}}\ln\sqrt{\frac{A'_nG'_n}{A_nG_n}}}<
\end{eqnarray}
$$
<\left(\frac{A_n}{G_n}\right)^{\frac{{A'_n}-{G'_n}}{A_n-G_n}-\frac{\ln\frac{A'_n}{G'_n}}
{\ln\frac{A_n}{G_n}}\ln\frac{A'_n}{A_n}}
<\left(\frac{A_n}{G_n}\right)^{1-\frac{\ln\frac{A'_n}{G'_n}}
{\ln\frac{A_n}{G_n}}\ln\frac{A'_n}{A_n}}<\frac{A_n}{G_n},
$$
and
\begin{eqnarray}
\frac{A'_n}{G'_n}<\max\left\{\left(\frac{A'_n}{G'_n}\right)^{1+\ln\sqrt{\frac{A'_nG'_n}{A_nG_n}}},
\left(\frac{A'_n}{G'_n}\right)^{\frac{A_n-G_n}{A'_n-G'_n}}\right\}<
\end{eqnarray}
$$
<\left(\frac{A'_n}{G'_n}\right)^{\left(1+\ln
\sqrt{\frac{A'_nG'_n}{A_nG_n}}\right)\frac{A_n-G_n}{A'_n-G'_n}}<\frac{A_n}{G_n},
$$
which are some refinements of Ky Fan's inequality {\rm (18)}.\\
Also,
\begin{eqnarray}
\frac{{A'_n}^n-{G'_n}^n}{{A_n}^n-{G_n}^n}<
\frac{\ln\left(\frac{A'_n}{G'_n}\right)^{{A'_n}^n{G'_n}^n}}{\ln\left(\frac{A_n}{G_n}\right)^{{A_n}^n{G_n}^n}},
\end{eqnarray}
which gives an inverse to {\rm (20)}.
\end{theo}
\begin{proof}
Since $A'_n>G'_n>A_n>G_n>0$, using the first inequality in (8), we
get the first inequality in (23) and the last one in (24). The
last inequality in (23), and using (19), the first one in (24) are
trivial. The other inequalities in (23) and (24) follow from (18)
and (19).\\ For proving (25), note that $f(-n)<f(0)$, where
$f(x)=\frac{{A'_n}^x-{G'_n}^x}{{A_n}^x-{G_n}^x}$.
\end{proof}
\begin{theo} Suppose $x_1, \cdots, x_n\in(0,\frac{1}{2}]$ not all equal. Then
\begin{eqnarray}
\max\left\{\frac{{A'_n}^n-{G'_n}^n}{{A_n}^n-{G_n}^n}\left(\frac{A_nG_n}{A'_nG'_n}\right)^{\frac{n}{2}},
\frac{A'_n-G'_n}{A_n-G_n}\left(\frac{A_nG_n}{A'_nG'_n}\right)^{\frac{1}{2}}
\right\}<
\end{eqnarray}
$$
<\frac{\ln\frac{A'_n}{G'_n}}{\ln\frac{A_n}{G_n}}<\frac{A'_n-G'_n}{A_n-G_n}
{\frac
{\ln\frac{I(A'_n,G'_n)}{I(A_n,G_n)}}{\ln\sqrt{\frac{A'_nG'_n}{A_nG_n}}}}<
\min\left\{\frac{A'_n-G'_n}{A_n-G_n},
{\frac{\ln\frac{I(A'_n,G'_n)}{I(A_n,G_n)}}{\ln\sqrt{\frac{A'_nG'_n}{A_nG_n}}}}
\right\}<1,
$$
and
\begin{eqnarray}
\frac{A'_n}{G'_n}<\left(\frac{A'_n}{G'_n}\right)^{\frac{\ln
\sqrt{\frac{A'_nG'_n}{A_nG_n}}}{\ln\frac{I(A'_n,G'_n)}{I(A_n,G_n)}}}
<\left(\frac{A_n}{G_n}\right)^{\frac{A'_n-G'_n}{A_n-G_n}}<\frac{A_n}{G_n}
<\left(\frac{A'_n}{G'_n}\right)^{\left(\frac{A'_nG'_n}{A_nG_n}\right)^{\frac{n}{2}}},
\end{eqnarray}
which give some refinements and inverses of Ky Fan's inequality {\rm (18)}.\\
Moreover,
\begin{eqnarray}
\frac{{A'_n}^n-{G'_n}^n}{{A_n}^n-{G_n}^n}<\frac{\ln(\frac{A'_n}{G'_n})^
{(A'_nG'_n)^{\frac{n}{2}}}}{\ln(\frac{A_n}{G_n})^{(A_nG_n)^{\frac{n}{2}}}}
<\left(\frac{A'_nG'_n}{A_nG_n}\right)^{\frac{n}{2}},
\end{eqnarray}
which gives some inverses of {\rm (20)}.
\end{theo}
\begin{proof} For proving the left hand side of (26), use
$\frac{G(a,b)}{G(c,d)}>\frac{L(a,b)}{L(c,d)}$ in Theorem 3.4 with
$a=A'_n,~ b=G'_n,~c=A_n,~d=G_n$ and also with $a={A'_n}^n,~
b={G'_n}^n,~c={A_n}^n,~d={G_n}^n$. For the second inequality in
(26), use (9) with $a=A'_n,~ b=G'_n,~c=A_n,~d=G_n$. The third and
forth inequalities in (26) follow from
$\frac{I(A'_n,G'_n)}{I(A_n,G_n)}<\sqrt{\frac{A'_nG'_n}{A_n,G_n}}$
in Theorem 3.4, and (19).\\
Putting $a=A'_n,~ b=G'_n,~c=A_n$ and $d=G_n$, the first inequality
in (27) follows from Theorem 3.4, the second one follows from (9),
the third one follows from (19),
and the last one follows from (26) by considering (20).\\
Inequalities in (28) follow from (26) and (18).
\end{proof}
\begin{theo} Suppose $x_1, \cdots, x_n\in(0,\frac{1}{2}]$ not all equal. Then
\begin{eqnarray}
\frac{A'_n}{G'_n}<\left(\frac{A_n}{G_n}\right)^{\frac{A_n+G_n}{A'_n+G'_n}\frac{A'_n-G'_n}{A_n-G_n}}<\min\left\{
\left(\frac{A_n}{G_n}\right)^{\frac{A_n+G_n}{A'_n+G'_n}},\left(\frac{A_n}{G_n}\right)^{\frac{A'_n-G'_n}{A_n-G_n}}\right\}<\frac{A_n}{G_n},
\end{eqnarray}
\begin{eqnarray}
\frac{A'_n}{G'_n}<\left(\frac{A'_n}{G'_n}\right)^{\frac{I(A'_n,G'_n)}{I(A_n,G_n)}\frac{A_n-G_n}{A'_n-G'_n}}
<\frac{A_n}{G_n},
\end{eqnarray}
which are some refinements of Ky Fan's inequality. Moreover, we
have
\begin{eqnarray}
\frac{A_nG_n}{A'_nG'_n}<\min\left\{
\left[\frac{\ln(\frac{A'_n}{G'_n})^
{\frac{1}{A'_n-G'_n}}}{\ln(\frac{A_n}{G_n})^{\frac{1}{A_n-G_n}}}\right]^2,
\left[\frac{\ln(\frac{{A'_n}^n}{{G'_n}^n})^
{\frac{1}{{A'_n}^n-{G'_n}^n}}}{\ln(\frac{{A_n}^n}{{G_n}^n})^{\frac{1}{{A_n}^n-{G_n}^n}}}\right]^{\frac{2}{n}}
\right\}<1.
\end{eqnarray}
\end{theo}
\begin{proof}
Put $a=A'_n,~ b=G'_n,~c=A_n,~d=G_n$. The inequalities in (29)
follow from $\frac{L(a,b)}{L(c,d)}>\frac{A(a,b)}{A(c,d)}$ in
Theorem 3.4, (19) and $A_n+G_n<A'_n+G'_n$.\\
The first inequality in (30) follows from Theorem 3.4 and (19),
and the second one follows from
$\frac{L(a,b)}{L(c,d)}>\frac{I(a,b)}{I(c,d)}$ in Theorem
3.4.\\
From (29) we obtain $(\frac{A'_n}{G'_n})^{\frac{1}{A'_n-G'_n}}
<(\frac{A_n}{G_n})^ {\frac{1}{A_n-G_n}}$, considering this
relation, we get (31) by solving the first inequality in (26) with
respect to $\frac{A_nG_n}{A'_nG'_n}$.
\end{proof}
\vspace{5mm} \textbf{Acknowledgment.} We would like to express our
gratitude to R. Kargar for his solution to the problem of
convexity of $f(x)=(a^x-b^x)/(c^x-d^x)$, which introduced us in
IASBS.

\end{document}